\documentclass[letter,12pt]{article}
\usepackage{amsfonts}
\usepackage{mathrsfs}
\usepackage{amsmath,amssymb,amsthm}
\usepackage{latexsym}
\usepackage{indentfirst}
\usepackage[top=1.25in, bottom=1.25in, left=1.25in, right=1.25in]{geometry}
\usepackage{color}

\let\footnote=\endnote

\newtheorem{cond}{Condition}
\newtheorem{rem}{Remark}
\newtheorem{thm}{Theorem}

\theoremstyle{definition}

\allowdisplaybreaks

\newcommand{{\m}}{{\cal \eta}}

\begin{document}

\title{ {\bf On Limiting Distribution of Quasi-Posteriors under Partial Identification}}
\author{ 
Wenxin Jiang\thanks{Taishan Scholar Overseas Distinguished Specialist Adjunct Professor at Shandong University, China, and Professor at Northwestern University,  U.S.A.; e-mail: wjiang@northwestern.edu.  
}
% cite yet without author's approval.
\rm\\
 Shandong University, China, and  Northwestern University,  U.S.A.
\rm}
\date{\today}
\maketitle
 \begin{center}
{\bf Abstract} \end{center} We establish the limiting distribution (in total variation) of the quasi posteriors based on moment conditions, which only   partially identify the parameters of interest. Some examples are discussed.   \\

\begin{small}
\noindent {\it Some key words:}   Generalized method of moments, interval data, moment inequalities,  partial identification, quasi-posterior, total variation.\\
\end{small}
MSC2010 Classification Codes: 62F15, 62F99.
\section{Introduction} 

Our paper studies theoretically the large sample behavior of certain Bayesian procedures which are of mutual interest to econometricians and statisticians. In Bayesian procedures, it is well known that the posterior distributions can often be approximated in total variation by normal distributions centered at the  frequentist maximum likelihood estimates (see, e.g., a popular textbook account of the Berstein-von Mises theorem in Chapter 10 of van der Vaart 2000). Econometricians have studied quasi-Bayesian approaches which require less assumptions by only assuming some moment conditions rather than a likelihood function (see, e.g., Kim 2002 and Chernozhukov and Hong 2003). In this context, similar and more general limiting results in, e.g., Chernozhukov and Hong (2003), suggest that the limiting posterior distribution is typically normal and centered at a corresponding frequentist extremum estimator such as one from the Generalized Method of Moments.

In the models with partial identification, where parameters are not point identified, so that the frequentist extremum estimator is not unique, the asymptotic normal limiting results mentioned before may fail. Such situations of partial identification have generated much interest recently both in statistics (e.g., Gustafson 2005, 2007, 2015) and in econometrics (e.g., Poirier 1998, Chernozhukov, Hong and Tamer 2007, Moon and Schorfheide 2012). The literature either focuses on the inference about the location of the partially identified   parameter (e.g., Moon and Schorfheide 2012, Gustafson 2015), or on the fully identified set of all possible locations of the parameter (e.g., Chernozhukov, Hong and Tamer 2007).  An incomplete sample of applications include missing data, interval censoring, (e.g., Manski and Tamer 2002, Manski 2003), game-theoretic models with multiple equilibria  (e.g., Bajari, Benkard, and Levin 2007, Ciliberto and Tamer 2009), auctions (Haile and Tamer 2003), noncompliance of randomized clinical trials (Gustafson 2015), and gene environment interactions (Gustafson 2015).
 
Our current paper derives and rigorously proves some results on the limiting posterior distribution in presence of partial identification.
In addition, we allow quasi-Bayes procedures based on moment conditions.
The limit is in the total variation sense, which is a Bernstein-von Mises type result, but not generally asymptotic normal.
When data are informative enough only to determine an identification region, instead of a point parameter, our result says that the limiting posterior is related to the prior distribution truncated in (a frequentist estimate of) this identification region. 

Our result  connects the literature on the inference about the identified parameter set to the inference of the unidentified parameter point, in the sense that one can easily convert a set estimate and combine it with a prior distribution to obtain a large sample approximation of the posterior distribution of the point parameter. 
This connection may have several meaningful applications.

1. (Simplifying computation)
 It can be used to avoid lengthy Markov chain Monte Carlo simulation that is typically involved in posterior computations, similar to using  the normal approximation in the point-identified situation. 

2. (Incorporating prior information.) This is also useful for incorporating prior information to improve the inference results from the more conservative set-based approach. The prior information can be from a same or different study with additional data that are either identifying or partially identifying the parameters of interest.   
For example, when a small part of the data are exact and all the rest are interval censored, 
it is obviously not advisable to use the interval data only to estimate an identification region. The exact part of the data can be used to derive a posterior, which can serve as a prior when further incorporating the interval data.

3. (Combining studies.)
In the above discussion, we have applied the principle that the prior can be derived from a posterior based on independent data. 
Multiple applications of this principle can allow meta-analysis (combining results from different studies), sequential computation (for dynamic data flow) or parallel computation (subsetting big data when they are hard to be handled altogether).  Our limiting posterior distribution suggests that combining inferences from subsets of data is equivalent to  intersecting their resulting identification regions.  

 \subsection{Related works.} 

A fundamental paper  about two decades ago by Poirier (1998) has shown many applications of the Bayesian method in handling problems with partial identification, where data $D$ are informative for only a subset of the parameters, say, $\lambda$, out of all the parameters $(\lambda,\theta)$ of the model. This decomposition of $(\lambda, \theta)$ does not have to be the most natural parametrization, but can be achieved by a clever re-parametrization. This situation is further illustrated by a sequence of works by Gustafson (e.g., 2005, 2007, 2015) with many interesting examples. The works of these authors have described that the limiting posterior distribution of $p(\lambda,\theta|D)=p(\lambda|D) p(\theta|\lambda)$ is the product of a usual asymptotic  normal distribution on $\lambda$, and a conditional prior $p(\theta|\lambda)$ where $\lambda$ either takes the true value or its maximum likelihood estimate. Recently this limiting result is rigorously proved in total variation distance by Moon and Schorfheide (2012, Theorem 1).
The key for this line of existing work is that (*) there exists a parameterization  decomposable into $(\lambda,\theta)$, such that the likelihood $p(D|\theta,\lambda)=p(D|\lambda)$ depends only on $\lambda$, i.e., $D$ and $\theta$ are independent given $\lambda$. Moon and Schorfheide (2012) call $\lambda$ the  ``reduced form'' parameter, and $\theta$ the ``structural'' parameter of interest.
 
 The current work aims for generalizing the works of these previous authors and studying the limiting posteriors under partial identification.
 The generalization is in two important ways:
 
{\em Generalization (i)(Posterior):} We generalize the likelihood-based posterior to be quasi-likelihood-based quasi-posterior, according to a general framework described in Chernozhukov and Hong (2003).

{\em Generalization (ii) (Partial identification):} We also allow more general scenarios of partial identification, where no obvious decomposition $(\lambda,\theta)$ can  satisfy (*), so that given $\lambda$, the data $D$ are $not$ ``conditionally'' uninformative / independent of the parameter of interest $\theta$.

To be more specific: we allow quasi-likelihood of the form $e^{-nR_n(\lambda,\theta)}$, where $R_n$ is a general empirical risk function that depends on data $D$, which was $(-1/n)$ times the log-likelihood function in the special case of the usual likelihood. We allow this quasi-likelihood to depend also on $\theta$ (unlike in (*) before), and only assume that there is a parametrization that can decompose into $(\lambda,\theta) $, such that (\dag) the marginal likelihood  $\int e^{-nR_n(\lambda,\theta)}d\lambda$ is a constant in $\theta$ (and is therefore ``marginally'' uninformative). It is obvious that this assumption (\dag) (marginal uninformativeness) contains (*)(conditional uninformativeness) as a special case, where $R_n=R_n(\lambda)$ had no dependence on $\theta$; but we will show later that there indeed exist interesting examples of (\dag) which do not satisfy (*).
 
\section{Examples}
%\bn
%\item Biased endogenous regression  example (likelihood or moment based)(more general than framework of reduced form parameter even in the likelihood case)
%\item Missing Data (likelihood based or moment based)
%\item Measurement Error
%\item Interval Censoring (mean model or quantile model)
%\item Inference of Classification Error (Liao and Jiang????)(marginal likelihood nontrivial due to nonconstant determinant)
%\item Collinearity
%\item Discrete classification (x is discrete and classification boundary is not point identified)(hope to cover)
%\item Identified case (hope to cover)
%\item cube asymptotics (not covered)
%%$\gb$
%\en

%TBA: 

%two perspectives of pid (cite Gustafson)(point vs set inference)

% compare with Moon and Schroefheide, and other literature review results.

%Cite Poirier and study his examples, Cite Gustafson and study his examples.
%cite Liao (2010)  thesis Theorem 3.4.2 (model and moment selection, not distributional limit in TV).

%Reduced form parameter special cases

%moment conditions

%first focus on tau=1.

%show their limit of mixture prior over posterior of reduced form parameter approaches our limit of conditional prior given true reduced form parameter.

%--

%WHAT ARE NEW

%more general partial identification setup then reduced form.

%more general quasi-likelihood

\subsection{Rounded data} 

This is based on a simple example of Section 2 of Moon and Schorfheide (2012). We here identify that it can be regarded as a special case of our framework. In this example, there is no Generalization (ii) in the structure of partial identification. Parameter decomposition (*) in the Introduction still holds. We use this example only to motivate the Generalization  (i) in using a quasi-likelihood (instead of a true likelihood).

Suppose we are interested in a structural parameter $\theta$ but it is only known that it is between
$[\phi,\phi+1]$, while we observe iid copies of $W\sim N(\phi,1)$.
Then the likelihood function is $p(data|\phi,\theta)\propto e^{-n R_n}\propto e^{-0.5\sum_{i=1}^n(W_i-\phi)^2}$, which is independent of the structural parameter $\theta$.  Here we can take the reduced parameter ($\lambda$ in the Introduction) to be $\phi$. 
Therefore (*) holds. Our  more relaxed condition (\dag) also holds: the marginalized likelihood 
$\int d\phi    e^{-nR_n}$ is totally uninformative of $\theta$, which fits our framework.

It is noted that we can use a quasi-likelihood based on moment conditions and be more flexible about the modeling of $W$. For example, suppose 
$\theta=EY$ where the observed $W=\lfloor Y \rfloor$ is the integer part of the hidden $Y$. Then $Y-W\in[0,1]$ and indeed $\theta\in [\phi,\phi+1]$ where $\phi=EW$. However, it is unsatisfactory to assume a normal model 
$W\sim N(\phi,1)$ since $W=\lfloor Y \rfloor$ is an integer.
In this case, our approach would only be to use the moment inequalities 
$ EW\leq EY\leq  EW+1$. We rewrite this with a moment equation 
$Em-\lambda=E(\theta-W)-\lambda=0$ where $\lambda$ $(=\theta-\phi)$ is a bias parameter constrained in $\lambda \in\Lambda=[0,1]$,
The corresponding sample moment is 
$\bar m-\lambda=\theta-\bar W-\lambda$.
Then we use a quasi-likelihood function
$e^{-n R_n}\propto e^{-0.5n (\bar m-\lambda)^2/v-0.5\log |2\pi v/n|} $ with $v=n^{-1}\sum_{i=1}^n(W_i-\bar W)^2$. 
The corresponding quasi-posterior will be
$q(\lambda,\theta|data)\propto e^{-n R_n} p(\lambda,\theta)I_{\Lambda\times\Theta}(\lambda,\theta)$, where $p(\lambda,\theta)$ is a prior and we have explicitly written out a region $\Lambda\times\Theta$ for   any constraints on $\lambda$ and $\theta$.
 The later development on the limiting distribution will show that  
 this quasi-posterior still makes sense in estimating the structural parameter $\theta$, even though we cannot easily use a genuine likelihood-based posterior in this case due to the rounded  values of the observed data $W$.
  
On the other hand, this simple example cannot be used to motivate Generalization (ii), since the structure of partial identification still satisfies the previous framework (*): we can  reparameterize and treat $\theta-\lambda$ $(=\phi)$ as the new reduced parameter $\lambda$. Given this new $\lambda$, data are conditionally uninformative about $\theta$. The next example shows that sometimes no obvious reparameterization like this exists, yet due to Generalization (ii), the example still satisfies the more general assumption (\dag) made in our proposed framework. 
 
%Therefore the marginalized likelihood (integrated without using the constraint $\Lambda$ from the prior information) is $\int_{\lambda\in \Re} d\lambda    e^{-nR_n}=1$ which is totally uninformative of $\theta$, which fits our framework.

\subsection{ Example of endogenous regression with biased error:}

This is a very messy example  which violates virtually all standard assumptions in linear regression.
Assume that the  observed response $Y$ in a   $nonlinear$ regression model follows $Y=g(X\theta)+\epsilon$ where $\epsilon$  is dependent on $X$ and has a marginal distribution $\epsilon\sim N(\lambda,\sigma^2)$    (with a scalar $X$ and a known $\sigma^2=1$ for simplicity), with a bias parameter $\lambda$ %(restricted to an interval $\Lambda$) 
describing the systematic component of the measurement error in $Y$.   (Later we will even allow nonnormality of the error.) Assume that there is no instrumental variable available, and we will consider inference based on the average $\bar Y$ over $n$ iid (independent and identically distributed) data: $\bar Y \sim N(\overline{g(X\theta)}+\lambda, 1/n)$. (Due to endogeneity, the likelihood function based on all data  is not guaranteed to be maximized near the true parameter. The sample average does not use the correlation with the endogenous variables, and therefore is still valid to use in inference.) 

This messy example  can be regarded as an extension of the 
basic case Example 3.1 of Poirier (1998), 
but it no longer leads to a identifiable reduced form parameter 
 ( such as $\theta+\lambda$  in the original example of a linear / constant mean model $Y=\theta+\epsilon$ with $\lambda=E\epsilon$ and  $g=\theta$), since the likelihood function based on $\bar Y$ is  proportional to $  e^{-0.5 n(\bar Y -\overline {g(X\theta)}-\lambda)^2}$, which is more complicated than before.
Given $\lambda$, data $\bar Y$  is $not$ conditionally independent of the structural parameter of interest $\theta$. So this example is now no longer covered by the framework of  Poirier (1998) or  Moon and Schorfheide (2012), even if we still used a   likelihood function to form the posterior.  On the other hand, this example is still covered by our framework.
%, suppose we use the likelihood 
%based on the average $\bar Y \sim N(\overline{g(X\theta)}+\lambda, 1/n)$, which is %proportional to
%$  e^{-0.5 n(\bar Y -\overline {g(X\theta)}-\lambda)^2}$.
This is because the marginalized likelihood
$ \int_{-\infty}^\infty d\lambda e^{-0.5 n(\bar Y -\overline {g(X\theta)}-\lambda)^2}=\sqrt{2\pi/n}$ is independent of the parameter of interest $\theta$.
In other words, Assumption (*) in Section 1 is not satisfied but Assumption (\dag) is, so we have an example of Generalization (ii) here.

Using the likelihood based on $\bar Y$ instead of based on the individual data $Y_i$'s has an additional advantage. The normal likelihood 
based on $\bar Y \sim N(\overline{g(X\theta)}+\lambda, 1/n)$ remains approximately valid for large $n$, even if the error $\epsilon_i$'s are $nonnormal$, due to the central limit theorem, as long as the lower order moments are correctly assumed.  In this case  it   really should be called a quasi-likelihood function.
Our framework can also cover such a more general situation where a quasi-likelihood function is  used to form a quasi-posterior (which was called Generalization (i) in Section 1).
%???NO!!!???.
 
\subsection{Interval Regression}
This follows from an extension of Example 2 of Chernozhukov, Hong and Tamer (2007). Assumes that $E(Y-g(X'\theta)|Z)=0$ for some positive instrumental variable $Z$ and a known parametric transform $g$. The structural parameter of interest is $\theta\in\Theta$. However, $Y$ is only observed to fall in an interval $[L,U]$.
This model cannot be easily treated in a reduce form. It is unclear how to form a reduced parameter so that data is independent of $\theta$ given the reduce parameter. In addition, it may not be desirable here to use a likelihood approach which would involve a joint probability model of $[L,U]$. Instead, we will use the moment conditions alone: $E(L|Z)\leq E(g(X'\theta)|Z)\leq E(U|Z)$, which implies $Em-\lambda\equiv (E(g(X'\theta)-L)Z^T-\lambda_1, E(U-L)Z-\lambda_2)^T=0$, constrained in 
$\lambda=(\lambda_1,\lambda_2)^T\in[0\leq \lambda_1\leq\lambda_2]\equiv \Lambda$.
 Consider a quasi posterior density of the form 
$$
q(\lambda,\theta|data)\propto e^{-nR_n(\lambda,\theta)}p(\lambda,\theta)I_\Xi(\lambda,\theta)
$$  
 where $p(\lambda,\theta)$ is a prior density on $\Xi=\Lambda\times\Theta$ with respect to a product base measure $d\lambda d\theta$, $\theta$ is the parameter of interest, and $\lambda$ is the nuisance parameter. Here 
$R_n$ is a  GMM (Generalized Method of Moments) criterion function $R_n=0.5  (\bar m(\theta)-\lambda)^Tv^{-1}(\bar m(\theta)-\lambda)+0.5n^{-1}\log |2\pi v/n|$ similar to the one used in Chapter 2 of Liao   (2010), where  $\bar m$ is the sample version of $Em$, and the variance matrix $v$ can be chosen as $I$ for simplicity (or alternatively as an estimate of $\sqrt{n}var \bar m$), which turns out to be irrelevant to the asymptotic inference about $\theta$ in the current partial identification scenario. 

Then this example is still covered in our framework, since it is obvious that the marginalized likelihood
$$
\int_{\lambda\in\Re^2}  e^{-nR_n(\lambda,\theta)} d\lambda =1,
$$ 
which is completely uninformative about $\theta$. In other words, the quasi-likelihood $e^{-nR_n(\lambda,\theta)}$ depends on $\theta$ and therefore violates (*) in Section 1, but (\dag) is still satisfied. This example therefore involves both Generalizations (i) and (ii).

%USE INTERVAL CENSORING INSTEAD??
 
\subsection{Interval Quantile Regression}\label{intrvl}

This example is similar to the previous one, except that we consider quantile regression here. Assumes that $E(q-I(Y\leq g(X'\theta))|Z)=0$ for some fixed $q\in(0,1)$, for some positive instrumental variable $Z$ and a known parametric transform $g$. The structural parameter of interest in $\theta\in\Theta$. However $Y$ is only observed to fall in an interval $[L,U]$.
This model cannot be easily treated in a reduce form. It is unclear how to form a reduced parameter so that data is independent of $\theta$ given the reduce parameter. In addition, it is not desirable here to use a likelihood approach which would involve a joint probability model of $[L,U]$. Instead we will use the moment conditions alone: $E(q-I(L\leq g(X'\theta))|Z)\leq 0$ and $E(q-I(U\leq g(X'\theta))|Z)\geq 0$.
 or alternatively, $0\leq E(q-I(U\leq g(X'\theta))|Z)  \leq  E(I(L\leq g(X'\theta))-I(U\leq g(X'\theta))|Z) $. This implies $Em-\lambda\equiv (E(q-I(U\leq g(X'\theta) )Z^T-\lambda_1, E(I(L\leq g(X'\theta))-I(U\leq g(X'\theta) ))Z-\lambda_2)^T=0$, constrained in 
$\lambda=(\lambda_1,\lambda_2)^T\in[0\leq \lambda_1\leq\lambda_2]\equiv \Lambda$.
 Consider a quasi posterior density of the form 
$$
q(\lambda,\theta|data)\propto e^{-nR_n(\lambda,\theta)}p(\lambda,\theta)I_\Xi(\lambda,\theta)
$$  
 where $p(\lambda,\theta)$ is a prior density on $\Xi=\Lambda\times\Theta$ with respect to a product base measure $d\lambda d\theta$, $\theta$ is the parameter of interest, and $\lambda$ is the nuisance parameter. Here 
$R_n$ is a GMM  criterion function $R_n=0.5  (\bar m(\theta)-\lambda)^Tv^{-1}(\bar m(\theta)-\lambda)+0.5n^{-1}\log |2\pi v/n|$ similar to the one used in  the previous example, with the new sample moment $\bar m$ corresponding to the quantile-related moment $Em$ defined above for the current example.

Then this example is still covered in our framework, since it is obvious that the marginalized likelihood
$$
\int_{\lambda\in\Re^2}  e^{-nR_n(\lambda,\theta)} d\lambda =1,
$$ 
which is completely uninformative about $\theta$. Therefore Assumption (\dag) in Section 1 is satisfied (but not (*)). This example therefore involves both Generalizations (i) and (ii). 
%USE INTERVAL CENSORING INSTEAD??

\subsection{Bayesian moment inequalities}\label{bmi}
Consider a quasi-posterior density of the form 
$$
q(\lambda,\theta|data)\propto e^{-nR_n(\lambda,\theta)}p(\lambda,\theta)I_\Xi(\lambda,\theta)
$$  
 where $p(\lambda,\theta)$ is a prior density on $\Xi=\Lambda\times\Theta$ with respect to a product base measure $d\lambda d\theta$, $\theta$ is the parameter of interest, and $\lambda$ is the nuisance parameter.
One example is the GMM   criterion function $R_n=0.5 (\bar m(\theta)-\lambda)^Tv^{-1}(\bar m(\theta)-\lambda)+0.5n^{-1}\log|2\pi v/n|$ similar to the one used in Chapter 2 of  Liao   (2010), where $v$ can be taken to be $I$ for simplicity (or alternatively  an estimate of $\sqrt{n}var \bar m$). When we constrain $\lambda\in\Lambda=[0,\infty)^{\dim(m)}$,
this corresponds to a Bayesian treatment of the moment inequality models which assume 
$Em(\theta)\geq 0$ (componentwise) (Chernozhukov, Hong and Tamer 2007), which can be rewritten as $Em(\theta)-\lambda=0$ subject to $\lambda\in\Lambda=[0,\infty)^{\dim(m)}$. The difference $\bar m(\theta)-\lambda$ used in $R_n$ is a sample version of $Em(\theta)-\lambda$. 
As pointed out Chernozhukov, Hong and Tamer (2007), the moment inequality models are useful for inference about the identification region 
$\{\theta\in\Theta: \ Em(\theta)\geq 0\ componentwise\}$ in many practically interesting examples.

The constraint $\lambda\in\Lambda$ may be regarded as part of the specification on the prior distribution $p(\lambda,\theta)I_{\Lambda\times\Theta} (\lambda,\theta)$. 
Before this prior constraint is imposed, the unconstrained integration over $\lambda$ of the quasi-likelihood $ e^{-nR_n(\lambda,\theta)}$ is uninformative about $\theta$, since
$$
\int_{\lambda\in\Re^{\dim(m)}}  e^{-nR_n(\lambda,\theta)} d\lambda =1.
$$
Therefore, (\dag) (but not (*)) in Section 1 is satisfied, and we have found another example   involving both Generalizations (i) and (ii).
 
%In Liao and Jiang (2010) Example 1.1, $\bar m$ is a moment function  
%$(n^{-1}\sum_1^n(\theta-Y_{1i}),n^{-1}\sum_1^n(Y_{2i}-\theta))$ and $\Lambda=[0,\infty)^2$.

\section{A general framework}

Consider a quasi posterior density of the form 
$$
q(\lambda,\theta|data)\propto e^{-nR_n(\lambda,\theta)}p(\lambda,\theta)I_\Xi(\lambda,\theta)
$$
where $p(\lambda,\theta)$ (when restricted to  $\Xi=\Lambda\times\Theta$) is a prior density  with respect to a product base measure $d\lambda d\theta$, $\theta$ is a parameter of interest, and $\lambda$ is a nuisance parameter. Here 
$R_n$ is an empirical risk function. 

The factor $e^{-nR_n(\lambda,\theta)}$ plays the role of a  likelihood function 
which summarizes the information from data, since the empirical risk $R_n$ depends on data.

(I) Assume that the likelihood function is $uninformative$ to the posterior inference on $\theta$, in the sense that the  {\em marginalized likelihood} $\int e^{-nR_n(\lambda,\theta)}d\lambda$, $before$ incorporating prior information, is constant in $\theta$, i.e., proportional to $\tau(\theta)=1$ (with an irrelevant proportional constant that can depend on $n$). Our examples  can all fit the choice $\tau(\theta)=1$.   More generally, we may allow the marginalized likelihood to be converging {\rm ``in some sense''} to a function proportional to  $\tau(\theta)$, which can $vary$ with $\theta$. Our theoretical results will be formally stated and proved in this more general framework. 

One possible way is to formalize this assumption as the following: 
\begin{cond}\label{condI}
There exists $C_n>0$ independent of $\lambda$ and $\theta$, and a nonstochastic function $\tau(\theta)$, such that
   $\tau_n(\theta)\equiv\int d\lambda C_n e^{-nR_n( \lambda ,\theta)}$ satisfies (\dag) 
$$
\int |\tau_n(\theta)/\tau(\theta)-1|^2p(\theta)d\theta=o_p(1).
$$
\end{cond}

(II) Assume that the   likelihood function is $informative$ to the posterior inference on $\lambda$ conditional on any given  $\theta$, in the sense that the usual Bayesian central limit theorem holds for $\lambda=\tilde\lambda(\theta)+t/\sqrt{n}$ around a first order extremum estimator $\tilde\lambda(\theta)$ (minimizing the empirical risk $R_n(\lambda,\theta)$ over $\lambda$ given $\theta$ asymptotically):
Given any $\theta$, the conditional density $$f(t|\theta)= e^{-nR_n(\lambda,\theta)-R_n(\tilde\lambda(\theta),\theta)}/\int \{e^{-nR_n(\lambda,\theta)-R_n(\tilde\lambda(\theta),\theta)}\} dt$$ converges in total variation and in probability to the normal  density of $N(0, V)$ where $V$ is a conditional asymptotic variance scaling as $1$ and can depend on $\theta$.  Some general conditions for this to happen are given by Belloni and Chernozhukov (2009, Theorem 1).

More formally, we assume the following (which may be provable under Belloni and Chernozhov's conditions for the posterior density under a flat prior, of $\lambda$ conditional on $\theta$, for almost all $\theta$ according to the prior $p(\theta)$):
\begin{cond}\label{condII} 
$$ \int dt  |  \{e^{-nR_n(\tilde\lambda(\theta)+t/\sqrt{n},\theta)}/\int ds e^{-nR_n(\tilde\lambda(\theta)+s/\sqrt{n},\theta)}  -\phi_V(t)\}   |=o_p(1),$$  for almost all $\theta$ according to the prior $p(\theta)$.  
\end{cond}

Under these two basic assumptions, and  with some additional mild regularity conditions, we have the following Theorem. 
\begin{thm} 
Under Conditions \ref{condI} and \ref{condII}, and additional mild regularity conditions \ref{posint},\ref{condidregion},\ref{avar}, \ref{condp}, \ref{bdry} (to be stated later), we have the following results:

(i) $p(\lambda,\theta|data)$ converges in total variation and in probability to a density proportional to 
$N(0, V/n)\times p(\tilde\lambda (\theta),\theta)I_{\Xi}(\tilde\lambda (\theta),\theta)\times \tau(\theta)$.

After integrating away the nuisance parameter $\lambda$, we
conjecture that 

(ii) $p(\theta|data)$  converges in total variation and in probability to a density proportional to 
$p(\tilde\lambda (\theta),\theta)I_{\Xi}(\tilde\lambda (\theta),\theta)\times \tau(\theta)$.

(iii) In the case   when the first order extreme estimator  $\tilde\lambda(\theta)$ 
converges   to a nonstochastic limit $\lambda(\theta)$ under Condition \ref{condlim},  
   the 
data dependent $\tilde\lambda(\theta)$ in $p(\tilde\lambda (\theta),\theta)I_{\Xi}(\tilde\lambda (\theta),\theta)$ that appears in both results (i) and  (ii) can be  replaced by $\lambda(\theta)$. More formally:
$$\int dt d\theta|e^{-nR_n(\lambda,\theta)}p(\lambda,\theta)I_\Xi(\lambda,\theta)/A-\phi_V(t) \tau(\theta)p(\lambda(\theta),\theta)I_\Xi(\lambda(\theta),\theta)/B|=o_p(1),$$
where $A=\int dt d\theta e^{-nR_n(\lambda,\theta)}p(\lambda,\theta)I_\Xi(\lambda,\theta)$
and $B=\int dt d\theta \phi_V(t) \tau(\theta)p(\lambda(\theta),\theta)I_\Xi(\lambda(\theta),\theta)$, and  $\lambda=\tilde\lambda(\theta)+t/\sqrt{n}$. The corresponding marginalized result is separately listed as:

(iv)
$$\int   d\theta|\int \lambda e^{-nR_n(\lambda,\theta)}p(\lambda,\theta)I_\Xi(\lambda,\theta)/ A- \tau(\theta)p(\lambda(\theta),\theta)I_\Xi(\lambda(\theta),\theta)/B|=o_p(1).$$
\end{thm}
Results (iii) and (iv) involve  a situation where the first order extremum estimator  $\tilde\lambda(\theta)$  converges in some sense to a nonstochastic limit $\lambda(\theta)$.  More formally,  they assume:
\begin{cond}\label{condlim}(On limit of the extremum estimator)
$$
 \int  d\theta   p(\theta) ||\tilde\lambda(\theta)-\lambda(\theta)||^2  = o_p(1).
$$
(Here $||\cdot||$ is the Euclidean norm.)
\end{cond}

Other mild regularity conditions include the following.

 \begin{cond}\label{condidregion}(On nondegenerate identification region.)
$$\int_\Theta p( \theta)I (\lambda(\theta)\in\Lambda) d\theta>0.$$
This means that the prior probability of an ``identification region'' $[\theta\in\Theta: \lambda(\theta)\in\Lambda]$ is positive.
\end{cond}

\begin{cond} \label{avar}(On regularity of the conditional assymptotic variance.) 
$$
\int  d\theta   p(\theta) |tr V|<\infty.
$$
\end{cond}

%\begin{cond} \label{avar}(On regularity of the conditional assymptotic variance.) 
%(\ddag)
%$$
%\int  d\theta   p(\theta) |tr V|/n  = o_p(1).
%$$
%\end{cond}

\begin{cond}\label{condp}(On regularity of the conditional prior.)
 (*)
$$
\int d\theta p(\theta) \tau(\theta)^2 [\sup_{\lambda } p(\lambda|\theta)^2]+\int d\theta p(\theta) \tau(\theta)^2[\sup_{\lambda } ||\partial_\lambda p(\lambda|\theta)||^2]<\infty.
$$
\end{cond}
\begin{cond}\label{bdry}(On prior probability of a boundary.) 
Let   $\lambda(\theta)$ be the large sample limit in  
Condition \ref{condlim}, and define a $\delta$-boundary $\partial_\delta  \Lambda=[\lambda: \max\{d(\lambda,\Lambda), d(\lambda,\Lambda^c)\}\leq\delta]$ for the region $ \Lambda $ (which is the parameter region of $\lambda$),
where $d$ between a point and a set is the minimal Euclidean distance. We assume that the prior distribution of $\lambda(\theta)$ is nonsingular on the boundary of $\Lambda$, i.e.,
(\$)
$$\lim_{\delta\downarrow 0}\int_\Theta d\theta  p(\theta)  I[\lambda(\theta) \in \partial_\delta  \Lambda]=0.$$ 
\end{cond}

\begin{rem}
\rm
Due to the Cauchy-Schwarz inequality, 
Conditions \ref{condidregion}  and \ref{condp} imply:
\begin{cond}\label{posint}(On positive normalizing constant of limiting density.)
$$\int p(\lambda(\theta),\theta)I_\Xi(\lambda(\theta),\theta)\tau(\theta)d\theta>0.$$
\end{cond}
Therefore, in results (iii) and (iv), the denominator $B>0$. [Also, $A>0$ with probability tending to 1, see this later in the proof of results (iii). Similarly in results (i) and (ii), the denominators in the normalizing constants are also positive with probability tending to 1.] 
\end{rem}

\begin{rem}
\rm The result (iii) can help explain why the limiting quasi-posterior makes sense relative to the true parameter relation $\lambda(\theta)$. The factor
$I_{\Lambda\times\Theta}(\lambda(\theta),\theta)$ can be recognized as the indicator function on the identification region of $\theta$. For example, in the moment inequality Example \ref{bmi}, one can easily verify that  $\tilde\lambda(\theta)=\bar m(\theta)$, $\lambda(\theta)=Em(\theta)$ and $\Lambda=[0,\infty)^{\dim(m)}$. Therefore $I_{\Lambda\times\Theta}(\lambda(\theta),\theta)=1$ only on the identification region $\{\theta\in\Theta: Em(\theta)\geq 0\; componentwise\}$. Results (iii) reasonably implies that   the limiting posterior gives 0 mass outside of the identification region and is prior dependent inside the identification region (and is proportional to $p(\theta)$ for flat $p(\lambda|\theta)$).
\end{rem}

\begin{rem}
\rm  The  result (i) and (ii) use data dependent $\tilde\lambda(\theta)$, which can be estimated by data (e.g., by $\bar m(\theta)$ in the moment inequality Example \ref{bmi}.) This has the advantage of obtaining a data-driven asymptotic distribution of the posterior. E.g., for flat $p(\lambda|\theta)$,
the posterior density  is asymptotically the same as  $constant\ p(\theta) I\{\theta: \bar m(\theta)
\geq 0\; componentwise\}$, which is the prior density truncated in an estimated identification region from a frequentist's approach. 
%(Discuss Gustafson 2014.) 
This can be used to compute the posterior distribution approximately without resorting to MCMC (Markov Chain Monte Carlo, as performed in Chapter 2 of Liao   2010).
% {example of linear instrumental regresson with interval censoring, polytope polyhedron???)
\end{rem}
\begin{rem}
\rm
 Both $\tau(\theta)$ and $V$ (which can depend on $\theta$ too) may be shown to be related to the second order derivatives of the large sample limit of $R_n$ in more general situations. However, there is no need to study these relations in detail in the current paper due to the following two reasons: (a). In all our examples in this paper, one can easily verify that  $\tau(\theta)$ can be simply taken to be 1. (b). The asymptotic variance  $V$ (for $\lambda$ conditional on $\theta$) does not affect  the limiting  posterior distribution marginally for $\theta$, which is often the only  parameter of interest. 
\end{rem}

%Remark: What is $\tau(\theta)$ under a quadratic approximation of $R_n(\lambda,\theta)$ given $\theta$? Discuss.
 
 %We need to do the following items of works:

%1. Find the conditions formally

%2. Prove the conjectures

%3. Find interesting examples

%4.  Verify the conditions for the examples

%5. Tell what useful things the general results imply for these examples. 

%6. Characterize the reasons and give examples why $\tau(\theta)$ can change with $\theta$:

%a. $R_\infty(\tilde\lambda(\theta),\theta)$ could change with $\theta$,

%b.  determinant $ |R{''}_\infty(\tilde\lambda(\theta),\theta)|^{-0.5}=V^{0.5}$ could change  with $\theta$,

%c. Jacobian when changing to more natural base measure from $d\lambda d\theta$ could depend on $\theta$ ???\\
%\newtheorem{thm2}{Proposition}
%\begin{thm2}
%\end{thm2}

\section{Proofs}
\subsection{Proof of results (iii).}
We first prove that for two nonnegative  functions $a,b$ such that $A=\int a>0$ and $B=\int b>0$ (with any common dominating measure suppressed in notation), we have 

\begin{equation}\label{ineq1}
\int |a/\int a-b/\int b|\leq 2\int |a-b|/ \int b.
\end{equation}

Proof of (\ref{ineq1}):

$\int |a/A-b/B| \leq \int (|b-a|/B+|a/B-a/A|)$

$= \int (|a-b|/B+|a||A-B|/(AB))$

$=\int|a-b|/B+|B-A|/B$
 
$\leq 2\int |a-b|/B$.

%We will prove (i).
%% for the case $\tau(\theta)=1$. 
%We will attempt to replace the data dependent $\tilde\lambda(\theta)$ in $p(\tilde\lambda (\theta),\theta)I_{\Xi}(\tilde\lambda (\theta),\theta)$ by the deterministic relation $\lambda(\theta)$ in the result, as formulated in the result (iii). [We will prove in the end that using the original formulation with $\tilde\lambda(\theta)$
%instead of $\lambda(\theta)$ is also OK.]
Then applying the  inequality (\ref{ineq1}) above, we can ignore the normalizing factor of $p(\lambda,\theta|data)$ and only need to prove that for some constant $C'_n>0$ independing of $\lambda$ and $\theta$,
$$
T\equiv \int dt d\theta|C'_ne^{-nR_n(\lambda,\theta)}p(\lambda,\theta)I_\Xi(\lambda,\theta)-\phi_V(t) \tau(\theta)p(\lambda(\theta),\theta)I_\Xi(\lambda(\theta),\theta)|=o_p(1).
$$
Here $\phi_V(t)=|2\pi V|^{-1/2}e^{-1/2t'V^{-1}t}$ is the density for $N(0,V)$,
and $\lambda=\tilde\lambda(\theta)+t/\sqrt{n}$. Note that Condition \ref{posint}  and $T=o_p(1)$ together imply that the two terms of the difference in $T$ both have positive integrals with probability tending to $1$ and therefore we can apply (\ref{ineq1}) to show that the normalized versions have difference $o_p(1)$.

We need to pay attention to the two indicator functions in this task, since the indicator function is not continuous in the usual sense. For this purpose,
now we introduce another inequality: If $I_{1,2}\in\{0,1\} $, then
\begin{equation}\label{ineq3}
\int |a I_1-b I_2| \leq \int |(a-b)I_1+b(I_1-I_2)|\leq \int |a-b| + \int |b||I_1-I_2|.
\end{equation}
 %> where $||b||_r=(\int |b|^r)^{1/r}$.  
Setting $I_{1,2}$ to be the two indicator functions in the results we wanted to prove above, we found that 
$$
T\leq T_1+T_2,
$$
and we only need to prove
$$
T_1\equiv \int dt d\theta|C'_ne^{-nR_n(\tilde\lambda(\theta)+t/\sqrt{n},\theta)}p(\lambda,\theta) -\tau(\theta)\phi_V(t) p(\lambda(\theta),\theta) |=o_p(1).
$$
and
$$
T_2\equiv \int dt d\theta\phi_V(t)\tau(\theta) p(\lambda(\theta),\theta)| I_\Xi(\lambda,\theta)- I_\Xi(\lambda(\theta),\theta)|=o_p(1).
$$
\subsubsection{Proof of $T_2=o_p(1)$.}

For $T_2$ with two indicator functions, we can rewrite it as 
$\int d\xi f(\xi)|I_\Lambda(\lambda)-I_\Lambda(\lambda(\theta))|=o_p(1)$,
where $\xi=(t,\theta)$, $\lambda=\tilde\lambda(\theta)+t/\sqrt{n}$,
 $f(\xi)=\phi_V(t) \tau(\theta)p(\lambda(\theta),\theta)I(\theta\in\Theta)$.

We split the integral domain into three parts: $A\cup B\cup  C$, 
where 

$A=[\xi: \max\{d(\lambda(\theta),\Lambda), d(\lambda(\theta),\Lambda^c)\}\leq\delta]$,

  $B=[\xi:  d(\lambda(\theta),\Lambda) >\delta]$,
  
  $C=[\xi:  d(\lambda(\theta),\Lambda^c)> \delta]$,
  
    for a minimal set distance under any metric $d$, and any $\delta>0$.   Notice that $I_\Lambda(\lambda(\theta))=0$ for $\xi\in B$ and $1$ for $\xi\in C$.

Then we bound the left hand side as follows: 
$\int d\xi f(\xi)|I_\Lambda(\lambda)-I_\Lambda(\lambda(\theta))|(I_A+I_B+I_C)
\leq \int d\xi f(\xi) I_A+\int d\xi f(\xi)I_\Lambda(\lambda) I_B+\int d\xi f(\xi)(1-I_\Lambda(\lambda))I_C$.
Now note that $I_\Lambda(\lambda)I_B$  and $(1-I_\Lambda(\lambda))I_C$ are indicator functions which can be one only when $d(\lambda,\lambda(\theta))>\delta$.
Then the integral is bounded by
$\int d\xi f(\xi)|I_\Lambda(\lambda)-I_\Lambda(\lambda(\theta))|\leq
\int d\xi f(\xi)  I_A+ 2\int d\xi f(\xi)I[d(\lambda,\lambda(\theta))>\delta]
\leq \int d\xi f(\xi)  I_A+ 2\int d\xi f(\xi) d(\lambda,\lambda(\theta)) /\delta $.
Therefore we have another inequality:
\begin{eqnarray}
&&\int d\xi f(\xi)|I_\Lambda(\lambda)-I_\Lambda(\lambda(\theta))| \leq \nonumber \\
&&\int d\xi f(\xi)  I[\max\{d(\lambda(\theta),\Lambda), d(\lambda(\theta),\Lambda^c)\}\leq\delta]+ 2\int d\xi f(\xi) d(\lambda,\lambda(\theta)) /\delta .\label{ineq4}
\end{eqnarray}
The first term in the upper bound will be related to the prior chance of $\lambda(\theta)$ falling within a distance of $\delta$ to the boundary of $\Lambda$, which does not depend on $n$, and is  typically converges to 0 as $\delta$ goes to 0. The second term is typically $O_p(1/\sqrt{n})/\delta$ (if $\tilde\lambda(\theta)$  is $\sqrt{ n}$-consistent for $\lambda(\theta)$).

The integrand of the first term in (\ref{ineq4}), after integrating away $t$, is equal to\\
$\int_\Theta d\theta \tau(\theta) p(\lambda(\theta),\theta)  I[\lambda(\theta) \in \partial_\delta  \Lambda]$ where   $\partial_\delta  \Lambda=A$.
This first integral satisfies \\ $\lim_{\delta\downarrow 0}\int_\Theta d\theta \tau(\theta)p(\lambda(\theta),\theta)  I[\lambda \in \partial_\delta  \lambda]=0$. 
   This is implied (after applying the Cauchy-Schwartz inequality) by   Condition \ref{condp} on boundedness  of the conditional prior, together with Condition \ref{bdry} which states that 
   
(\$)
$ \lim_{\delta\downarrow 0}\int_\Theta d\theta  p(\theta)  I[\lambda(\theta)  \in \partial_\delta  \Lambda]=0.$

For the second term in (\ref{ineq4}), choose $d$ to be the Euclidean metric $||\cdot||$. Then 
 $d(\lambda,\lambda(\theta)) =  ||\tilde\lambda(\theta)+t/\sqrt{n}-\lambda(\theta)|| \leq      (||\tilde\lambda(\theta) -\lambda(\theta)||  +||t/\sqrt{n}|| ) $. Then
the second term in (\ref{ineq4}) is bounded by 

$2 \int_\Theta d\theta \tau(\theta)p(\lambda(\theta),\theta) [||\tilde\lambda(\theta)-\lambda(\theta)||  + \int ||t|| \phi_V(t)dt/\sqrt{n}]/\delta $

 $ \leq \sqrt{\int d\theta p(\theta) \tau(\theta)^2p(\lambda(\theta)|\theta)^2}[\sqrt{  \int d\theta p(\theta)||\tilde\lambda(\theta)-\lambda(\theta)||^2 } +\sqrt{  \int d\theta p(\theta)  |tr V|/n}]/\delta = o_p(1)/\delta$.
[This is implied by Conditions \ref{condp} and \ref{condlim}.
We can  also apply Condition \ref{avar} regarding the $V$ matrix.] 
 
% We will assume:
%\begin{enumerate}
%\item Metric $d$ is  equivalent to the Euclidean metric $||\cdot||$.
%\item The first integral satisfies $\lim_{\delta\downarrow 0}\int d\theta \tau(\theta)p(\m)  I[\m \in \partial_\delta  \Xi]=0$,
% where $\m=(\lambda(\theta),\theta)$ and $\partial_\delta  \Xi=[\theta: \max\{d(\m,\Xi), d(\m,\Xi^c)\}\leq\delta]$.???
%[This will be implied by a boundedness condition (*)   of the conditional prior  appearing later, together with (\$)
%$$\lim_{\delta\downarrow 0}\int d\theta  p(\theta)  I[(\lambda(\theta),\theta) \in \partial_\delta  \Xi]=0.$$]
%\item  The integral $\int d\theta \tau(\theta)p(\m) ||\tilde\lambda(\theta)-\lambda(\theta)||^2=o_p(1)$. [This will be implied by (*) and (**) appearing later.]
%\item A later assumption (\ddag) regarding $V$ matrix.
%\end{enumerate}

Then we have proved that 
$\int dt d\theta\phi_V(t) \tau(\theta) p(\lambda(\theta),\theta)| I_\Xi(\lambda,\theta)- I_\Xi(\lambda(\theta),\theta)|=o_\delta(1)+o_p(1)/\delta$  for any small positive $\delta$,
where $\lim_{\delta\downarrow 0}o_\delta(1)=0$. Therefore
$$
 T_2\equiv\int dt d\theta\phi_V(t) \tau(\theta) p(\lambda(\theta),\theta)| I_\Xi(\lambda,\theta)- I_\Xi(\lambda(\theta),\theta)|=o_p(1).
$$
  
\subsubsection{Proof of $T_1=o_p(1)$.}

Now we return to  the other statement involving $T_1$ without the indicator functions. We wanted to prove that there exists $C'_n>0$ independent of $\lambda$ and $\theta$, such that
$$
T_2\equiv \int dt d\theta|C'_ne^{-nR_n(\tilde\lambda(\theta)+t/\sqrt{n},\theta)}p(\lambda,\theta) -\phi_V(t) \tau(\theta)p(\lambda(\theta),\theta) |=o_p(1).
$$
(Unless otherwise noted, $\lambda$ is related to $t$ by the reparameterization $\lambda=\tilde \lambda(\theta)+t/\sqrt{n}$.)
%The left hand side is to be re-written as\\
%$
%\int dt d\theta|\{\int ds e^{-nR_n(\tilde\lambda(\theta)+s/\sqrt{n},\theta)} \} \{e^{-nR_n(\tilde\lambda(\theta)+t/\sqrt{n},\theta)}/\int ds e^{-nR_n(\tilde\lambda(\theta)+s/\sqrt{n},\theta)} \}p(\lambda,\theta) -\phi_V(t) p(\lambda(\theta),\theta) | =\int dt d\theta|  \{e^{-nR_n(\tilde\lambda(\theta)+t/\sqrt{n},\theta)}/\int ds e^{-nR_n(\tilde\lambda(\theta)+s/\sqrt{n},\theta)} \}p(\lambda,\theta) -\phi_V(t) p(\lambda(\theta),\theta) |
%.
%$
%This last step is obtained from the assumption that 
%the  factor $$\int ds e^{-nR_n(\tilde\lambda(\theta)+s/\sqrt{n},\theta)}=\tau(\theta)=1.$$

%More generally, 

To rewrite the left hand side, we can let $C'_n=C_n/n^{\dim(\lambda)/2}$ for $C_n$ in Condition \ref{condI}, then 
 $\int ds C'_n e^{-nR_n(\tilde\lambda(\theta)+s/\sqrt{n},\theta)}=\tau_n(\theta) $, and
%n^{\dim(\lambda)/2}

 $f(t|\theta)=e^{-nR_n(\tilde\lambda(\theta)+t/\sqrt{n},\theta)}/\int ds e^{-nR_n(\tilde\lambda(\theta)+s/\sqrt{n},\theta)} $. Then 
the left hand side can be re-written as\\
$T_1=$

$
\int dt d\theta|[\{\tau_n(\theta)-\tau(\theta) \} +\tau(\theta)]f(t|\theta)p(\lambda,\theta) -\tau(\theta)\phi_V(t) p(\lambda(\theta),\theta) |$

$\leq \int |\tau_n(\theta)-\tau(\theta)|p(\theta)\sup_{\lambda }|p(\lambda|\theta)|d\theta +  \int dt d\theta \tau(\theta) |  f(t|\theta)p(\lambda,\theta) -\phi_V(t) p(\lambda(\theta),\theta) |$

$\equiv T_{11}+T_{12}.$

The first term $T_{11}=o_p(1)$ due to boundedness of the conditional prior density from Condition \ref{condp} (*), and the relation (\dag) in Condition \ref{condI} that
$$
\int |\tau_n(\theta)/\tau(\theta)-1|^2p(\theta)d\theta=o_p(1).
$$
Now we rewrite

$T_{12}\equiv
\int dt d\theta \tau(\theta)|  \{e^{-nR_n(\tilde\lambda(\theta)+t/\sqrt{n},\theta)}/\int ds e^{-nR_n(\tilde\lambda(\theta)+s/\sqrt{n},\theta)} \}p(\lambda,\theta) -\phi_V(t) p(\lambda(\theta),\theta) |$

$\leq \int dt d\theta|  \tau(\theta) \{e^{-nR_n(\tilde\lambda(\theta)+t/\sqrt{n},\theta)}/\int ds e^{-nR_n(\tilde\lambda(\theta)+s/\sqrt{n},\theta)}  -\phi_V(t)\} p(\lambda,\theta)  |$

$
+\int dt d\theta \tau(\theta)\phi_V(t)|(p(\lambda,\theta)-p(\lambda(\theta),\theta)) |$

$\equiv T_{121}+T_{122}$,

where $\lambda=\tilde \lambda (\theta)+t/\sqrt{n}$.

Now assume the following ((*) from Condition \ref{condp}):
$$
\int d\theta p(\theta) \tau(\theta)^2 [\sup_{\lambda } p(\lambda|\theta)^2]+\int d\theta p(\theta) \tau(\theta)^2[\sup_{\lambda } ||\partial_\lambda p(\lambda|\theta)||^2]<\infty.
$$
Then the second term

 $T_{122}=$ 

$\int dt d\theta \tau(\theta)\phi_V(t)|(p(\lambda,\theta)-p(\lambda(\theta),\theta)) |$

$
\leq \int dt d\theta \phi_V(t)\tau(\theta)p(\theta)[\sup_{\lambda } ||\partial_\lambda p(\lambda|\theta)||] ||\lambda-\lambda(\theta)||$ 

$
\leq  \int dt d\theta \phi_V(t) \tau(\theta))[\sup_{\lambda } ||\partial_\lambda p(\lambda|\theta)||] p(\theta)( ||\tilde\lambda(\theta)-\lambda(\theta)||+||t||/\sqrt{n})  
$

$\leq   
 \sqrt{\int   d\theta  p(\theta)( ||\tilde\lambda(\theta)-\lambda(\theta)||^2  +|tr V|/n)} \sqrt{\int   d\theta  p(\theta) \tau(\theta)^2[\sup_{\lambda } ||\partial_\lambda p(\lambda|\theta)||] ^2}$
 
 $
= 
o_p(1),$

assuming    (**) (from Condition \ref{condlim})

$$
 \int  d\theta   p(\theta) ||\tilde\lambda(\theta)-\lambda(\theta)||^2  = o_p(1),
$$
and (from Condition \ref{avar}) noting that
(\ddag) 
$$
\int  d\theta   p(\theta) |tr V|/n  = o_p(1).
$$

For the first term,

$T_{121}=  \int dt d\theta|  \{e^{-nR_n(\tilde\lambda(\theta)+t/\sqrt{n},\theta)}/\int ds e^{-nR_n(\tilde\lambda(\theta)+s/\sqrt{n},\theta)}  -\phi_V(t)\}\tau(\theta)p(\lambda,\theta)  |$

$\leq 
 \int dt d\theta|  \{e^{-nR_n(\tilde\lambda(\theta)+t/\sqrt{n},\theta)}/\int ds e^{-nR_n(\tilde\lambda(\theta)+s/\sqrt{n},\theta)}  -\phi_V(t)\}| p(\theta) \tau(\theta) \sup_{ \lambda } p(\lambda|\theta)  \equiv B_1.$

Now we assume Condition \ref{condII},  which states that
$$ \int dt  |  \{e^{-nR_n(\tilde\lambda(\theta)+t/\sqrt{n},\theta)}/\int ds e^{-nR_n(\tilde\lambda(\theta)+s/\sqrt{n},\theta)}  -\phi_V(t)\}   |=o_p(1),$$  for almost all $\theta$ according to the prior $p(\theta)$.    Then for all these $\theta$, 
$$
E \int dt  |  \{e^{-nR_n(\tilde\lambda(\theta)+t/\sqrt{n},\theta)}/\int ds e^{-nR_n(\tilde\lambda(\theta)+s/\sqrt{n},\theta)}  -\phi_V(t)\}   |=o(1)
$$
since the $L_1$ distance of two densities is bounded by 2, and convergence in probability implies convergence in mean.
 Then by using a dominated convergence theorem and noting that $E\int dt  |  e^{-nR_n(\tilde\lambda(\theta)+t/\sqrt{n},\theta)}/\int ds e^{-nR_n(\tilde\lambda(\theta)+s/\sqrt{n},\theta)}  -\phi_V(t) |\leq 2$, which is integrable under $\int (...) p(\theta)\tau(\theta) \sup_{ \lambda } p(\lambda|\theta) d\theta$ due to assumption (*) from Condtion \ref{condp}, we arrive at
$$
\int (E \int dt  |  \{e^{-nR_n(\tilde\lambda(\theta)+t/\sqrt{n},\theta)}/\int ds e^{-nR_n(\tilde\lambda(\theta)+s/\sqrt{n},\theta)}  -\phi_V(t)\}   |) p(\theta)\tau(\theta) \sup_{ \lambda } p(\lambda|\theta) d\theta =o(1).
$$  
Exchanging $\int d\theta$ and $E$ by Fubini's theorem, we obtain

$
 E \int (\int dt  |  \{e^{-nR_n(\tilde\lambda(\theta)+t/\sqrt{n},\theta)}/\int ds e^{-nR_n(\tilde\lambda(\theta)+s/\sqrt{n},\theta)}  -\phi_V(t)\}   |) p(\theta)\tau(\theta) \sup_{ \lambda } p(\lambda|\theta) d\theta 
$
$=
EB_1=o(1).
$
Then $0\leq T_{121}\leq B_1=o_p(1)$.

So far, the arguments above, when collected together, have proven that $T_1=o_p(1)$. This, together with the results from earlier subsections, have proven the theorem in the formulation of result (iii), using the deterministic relation $\lambda(\theta)$ in the result.  Q.E.D.

\subsection{Proof of results (i).}

The original result  (i) is formulated with  the data dependent $\tilde\lambda(\theta)$ in $p(\tilde\lambda (\theta),\theta)I_{\Xi}(\tilde\lambda (\theta),\theta)$. We will now prove   that using   $\tilde\lambda(\theta)$
instead of $\lambda(\theta)$ is also OK, in the sense that the limiting densities differ only by $o_p(1)$ in the $L_1$-distance.  Due to (\ref{ineq1}), it suffices for us to prove that 
$$
S\equiv \int dt d\theta|\phi_V(t) p(\tilde\lambda(\theta),\theta)I_\Xi(\tilde \lambda(\theta),\theta)\tau(\theta)-\phi_V(t) p(\lambda(\theta),\theta)I_\Xi(\lambda(\theta),\theta)\tau(\theta)|=o_p(1).
$$
The variable $t$ can be integrated away.

%\begin{cond}\label{posint}
%$\int p(\lambda(\theta),\theta)I_\Xi(\lambda(\theta),\theta)\tau(\theta)d\theta>0$.
%\end{cond}

By Condition \ref{posint} $\int p(\lambda(\theta),\theta)I_\Xi(\lambda(\theta),\theta)\tau(\theta)d\theta>0$, and $S=o_p(1)$, both terms in the difference of $S$ should have positive integrals with probability tending to 1, which enables us to use (\ref{ineq1}) and show that the normalized versions also have difference $o_p(1)$.

Now we use again the aforementioned
  inequality (\ref{ineq3}): If $I_{1,2}\in\{0,1\} $, then
$\int |a I_1-b I_2| \leq \int |a-b| + \int |b||I_1-I_2|$. %> where $||b||_r=(\int |b|^r)^{1/r}$.  
Setting $I_{1,2}$ to be the two indicator functions in the results we wanted to prove above, we found that 
$$
S\leq S_1+S_2,
$$
and we only need to prove
$$
S_1\equiv  \int d\theta \tau(\theta)|p(\tilde\lambda(\theta),\theta)-  p(\lambda(\theta),\theta) |=o_p(1),
$$
and
$$
S_2\equiv \int d\theta  \tau(\theta)p(\lambda(\theta),\theta)| I_\Xi(\tilde\lambda(\theta),\theta)- I_\Xi(\lambda(\theta),\theta)|=o_p(1).
$$
 
For the first term, 
$$
S_1=\int  d\theta \tau(\theta)|(p(\tilde \lambda(\theta),\theta)-p(\lambda(\theta),\theta)) |
\leq \int d\theta  p(\theta)\tau(\theta) [\sup_{\lambda } ||\partial_\lambda p(\lambda|\theta)||]||\tilde\lambda(\theta)-\lambda(\theta)||,$$ 
which is $ o_p(1)$ due to the assumptions   made before from Conditions  \ref{condp} and \ref{condlim}.

For the second term $S_2$ we apply again the aforementioned  inequality (\ref{ineq4}):
\begin{eqnarray}\label{eq.temp}
&&\int d\xi f(\xi)|I_\Lambda(\lambda)-I_\Lambda(\lambda(\theta))| \leq \nonumber \\
&&\int d\xi f(\xi)  I[\max\{d(\lambda(\theta),\Lambda), d(\lambda(\theta),\Lambda^c)\}\leq\delta]+ 2\int d\xi f(\xi) d(\lambda,\lambda(\theta)) /\delta . 
\end{eqnarray}
%
%$$
%\int d\xi f(\xi)|I_\Xi(\hat \m(\xi))-I_\Xi(\m(\xi))| \leq$$
%$$\int d\xi f(\xi)  I[\max\{d(\m(\xi),\Xi), d(\m(\xi),\Xi^c)\}\leq\delta]+ 2\int d\xi f(\xi) %d(\hat \m(\xi),\m(\xi))^2/\delta^2.
%$$
Now we take $\xi= (\theta,t)$, $\lambda=\tilde\lambda(\theta)$,  
 $f(\xi)=  p(\lambda(\theta),\theta)\tau(\theta)I(\theta\in\Theta)\phi_V(t)$.
Then the left hand side of (\ref{eq.temp}) can be recognized to be $S_2$, and its upper bound is $o_\delta(1)+o_p(1)/\delta$  due to Conditions \ref{condlim}, \ref{bdry} and \ref{condp}, where $o_\delta(1)$ converges to 0 as $\delta\downarrow 0$, and is independent of data.  This shows that
the second term $S_2$ is also $o_p(1)$. 

Collecting the arguments above, we have shown that result (i) also holds with the data dependent relation $\tilde\lambda(\theta)$ used (instead of the deterministic $\lambda(\theta)$) in the result. 
 Q.E.D.
 \subsection{Proof of results (ii) and (iv).}

The total variation (or $L_1$) distance of the joint densities is stronger than that of the corresponding marginal densities, i.e., 
\begin{equation}\label{ineq2}
\int|q(\lambda,\theta)-q'(\lambda,\theta)|d\lambda d\theta \geq \int|\int q(\lambda,\theta)d\lambda-\int q'(\lambda,\theta)d\lambda| d\theta .\end{equation}
 Therefore,  the convergence of the joint distributions (result (i) and result (iii)) implies  the convergence of the marginal distributions (results  (ii) and (iv), respectively). Q.E.D.

%(TBA: list all conditions separately before hand. Formulate Lemmas if needed.)
%\\
\section{Regularity conditions for BGMM (Bayesian Generalized Method of Moments)}
We now study the 7 regularity conditions for a general class of quasi-posteriors obtained from BGMM. 
Consider a quasi posterior density of the form 
$$
q(\lambda,\theta|data)\propto e^{-nR_n(\lambda,\theta)}p(\lambda,\theta)I_\Xi(\lambda,\theta)
$$  
 where 
$R_n$ is an  GMM (Generalized Method of moments) criterion function $$R_n=0.5  (\bar m(\theta)-\lambda)^Tv^{-1}(\bar m(\theta)-\lambda)+0.5n^{-1}\log |2\pi v/n|,$$ similar to the one used in  Chapter 2 of Liao   (2010), where  $\bar m$ is the sample version of $Em$, and the variance matrix $v$ can be chosen as $I$ for simplicity (or alternatively by an estimate of $\sqrt{n}var \bar m$), which turns out to be irrelevant to the asymptotic inference about $\theta$ in the current partial identification scenario.  When we need a more explicit form, we will consider a sample average $\bar m(\theta)=n^{-1}\sum_{i=1}^n m(W_i,\theta)$, and $Em(\theta)=Em(W,\theta)$, where $W, W_1,...,W_n,...$ are iid (independent and identically distributed).

Condition \ref{condI}: It is obvious that the marginalized likelihood
$$
\int   e^{-nR_n(\lambda,\theta)} d\lambda =1.
$$ 
 So Condition \ref{condI} is satisfied with $C_n=1$ and $\tau(\theta)=1$.

 Condition \ref{condII}: the extremum estimator of $\lambda$ given $\theta$ is $\tilde\lambda(\theta)=\bar m(\theta)$. With the BGMM choice of $R_n$, the quasi-likelihood  $e^{-nR_n}$ is already proportional to a normal density $\phi_v$ with variance $v$.
 We only need that $\phi_v$ converges to $\phi_V$ in total variation. This is achievable  if $v$ is consistent estimator of $V$. This happens, when, e.g., $v$ is a sample version of $V= var_{W|\theta}   m(W, \theta)$.

 Condition \ref{condlim}: We can take $\lambda(\theta)=Em(\theta)$. Then the condition is satisfied when $p(\theta)$ is supported on a compact set $\Theta$, and when $\bar m(\theta)$ converges to $Em(\theta)$ uniformly on $\Theta$, in probability. 
 
Condition \ref{condidregion}:
This condition means that the prior probability of the identification region $[\theta\in\Theta: \lambda(\theta)\in\Lambda]$ is positive. This is reasonable in many partial identification problems where the identification region is nondegenerate.

 Condition \ref{avar}:
     When $tr(V)$ is bounded in $\theta$ on the support of $p(\theta)$, it is obvious that the integral is finite.
     
 Condition \ref{condp}:              
We do not need to worry about the $\tau(\theta)$ factor, since we have $\tau(\theta)=1$ for BGMM. 
Suppose there exists an extension of $p(\lambda|\theta)$   from $\Lambda\times\Theta$ to  $\Re^{\dim(\lambda)}\times\Theta$ such that its function values and $\partial_\lambda$ derivatives are all bounded functions, then the condition is obviously satisfied.

     Condition \ref{bdry}:
     Note that $\partial_\delta\Lambda$ typically has Lesbegue measure $O(\delta)$ in the direction of one $\lambda_j$ component.  For Example \ref{bmi} (with the Bayesian moment inequalities), $\Lambda=[0,\infty)^{\dim(m)}$.  The event $\lambda(\theta)\in \partial_\delta\Lambda$ implies that some $|\lambda_j(\theta)|\leq \delta$ for some $j\in\{1,...,\dim(m)\}$. As long as the prior density  for $\lambda_j(\theta)$ is finite at $\lambda_j(\theta)=0$ for all $j$, its integral on $[-\delta,+\delta]$ will be $O(\delta)$, guaranteeing that the condition holds.
     
     By basic calculus, the prior density of $\lambda_j(\theta)$ can be computed by reparameterizing it in a form of  $\int  d\theta_{(-k)}p(\theta_k(\lambda_j, \theta_{(-k)}),\theta_{(-k)})/|\partial_{\theta_k}\lambda_j(\theta)|$,
     for some decomposition of $\theta$ into some $\theta_k$ and all other components $\theta_{(-k)}$. Suppose that the prior density $p(\theta)$ is bounded and supported on a 
  bounded set $\Theta$, and that $   |\partial_{\theta_k}\lambda_j(\theta)|$ is bounded away from 0 on $\theta\in\Theta$. Then the prior density of $\lambda_j(\theta)$ is finite.
  
 This last condition on the derivative can be verified by noting that $\lambda(\theta)=Em(\theta)$. For example, suppose $Em=[EZ(e^{X'\theta}-L), EZ(U-e^{X'\theta})]$ for some positive instrumental variable  $Z$, as is useful for a  Bayesian moment inequality approach of the interval regression model $E(Y-e^{X'\theta}|Z)=0$ where $Y$ is only known to fall in $[L,U]$.
 Then the absolute value of the derivative of any component $Em$ against one chosen $\theta$ component is of the form $ | EZ_kX_ke^{X'\theta}|$. Suppose $Z_kX_k>0$ (e.g., suppose we can translate $X_k$ to make $X_K>0$ and we take $Z_k=X_k$), then $| EZ_kX_ke^{X'\theta}|\geq |E|Z_kX_k|e^{-\sup|X'\theta|} >0$, if we assume bounded $X'\theta$.  Then the derivative condition (and therefore Condition \ref{bdry}) is satisfied.
 
 %It can be verified in all examples given in this paper. We just describe one example since other examples are similar. For   Example \ref{intrvl} (with interval regression),   $\lambda(\theta)=Em(\theta)$, which has a component 
% . The event $\lambda(\theta)\in \partial_\delta\Lambda$ implies that some $|\lambda_j(\theta)|\leq \delta$
                 
% $\lambda(\theta)=Em(\theta)$.
\section{Discussions}
In this paper, we have derived the limiting distribution (in total variation) of the posterior distribution under partial identification. Our proof is rigorous, and the framework is general enough to include quasi-Bayes methods based on moment conditions. In addition, we allow more general partial identification, where the model may not be easily reparameterized to be an identifiable model with some reduced form parameters. The resulting limit of the posterior distribution combines information from the data and from the prior  reasonably: it uses the data information only to locate an identifiable region, and then leaves the within-region knowledge to be determined by the prior distribution.
 
In the Bayesian literature of partial identification,  there is a new  direction of work where the Bayesian inference is targeted at the identification set, rather than a point parameter. See, for example, Kline and Tamer (2016), and Chen, Christensen and Tamer (2016). This direction is different and interesting, and has the advantage of stating conclusions robustly without being influenced by additional assumptions   on the prior distributions. Another earlier work (Kitagawa  2012) explicitly addresses this robustness aspect associated with targeting at the identification set, using bounds on the posterior probabilities due to a class of priors.  Our current paper, on the other hand,  uses the traditional framework of Bayesian inference, in the sense that the unknown true parameter is regarded as a random point in a parameter space. This follows the line of work by Poirier (1998), Gustafson (2005, 2007, 2015)), and Moon and Schorfheide (2012), and has the advantage of being able to improve the parametric inference by incorporating useful prior information.  
Both approaches are content on accepting partial identification and  
are robust regarding the mechanism of missing data, as compared to other approaches that strive for point identification by introducing additional assumptions on the missing data mechanism.

\vspace*{0.5cm}

{\bf Acknowledgments}

I thank Professor Hyungsik Roger Moon for kindly reading a draft of this paper and providing useful references.

% {\bf REFERENCES TO BE TIDIED UP and CHECKED.}


\begin{thebibliography}{99}

\bibitem{} Bajari, P., L. Benkard, \& J. Levin (2007). Estimating dynamic models of imperfect competition, {\it Econometrica} 75, 1331–1370.  


\bibitem{} Belloni, A. \& V. Chernozhukov (2009). On the computational complexity of MCMC-based estimators in large
samples. \textit{The Annals of Statistics} 37, 2011-2055.

%\bibitem{} Bochkina, N. A. \& P. J. Green (2014). The Bernstein-von Mises theorem and nonregular
%models. The Annals of Statistics 42(5), 1850–1878.

\bibitem{} Ciliberto, F., \& E. Tamer (2009). Market structure and multiple equilibria in airline markets, {\it Econometrica} 77, 1791–1828.  


\bibitem{} Chen, X. Christensen, T. \&  E. Tamer (2016).
MCMC confidence sets for identified sets.
{\it Cowles Foundation Discussion Paper} No. 2037.
http://papers.ssrn.com/sol3/papers.cfm?abstract\_id=2775253 

\bibitem{} Chernozhukov, V. \&  H. Hong (2003). An MCMC approach to classical estimation.
\textit{Journal of Econometrics} 115, 293-346.

\bibitem{} Chernozhukov, V.,   Hong, H.  \& E. Tamer  (2007). Estimation and conﬁdence regions for parameter sets in econometric models. {\it Econometrica} 75, 1243–1284. 


\bibitem{}Gustafson, P. (2005). On model expansion, model contraction, identifiability, and prior information: two illustrative
scenarios involving mismeasured variables (with discussion). {\it Statist. Sci.}, 20, 111–140.

\bibitem{}Gustafson, P. (2007).
Measurement error modelling with an approximate
instrumental variable. {\it J. R. Statist. Soc. B} 69,  797–815.
 
% \bibitem{}Gustafson, P. (2014). Bayesian inference in partially identified models: Is the shape of the posterior distribution useful?
 %{\it Electronic Journal of Statistics}  8, 476–496.

\bibitem{}Gustafson, P. (2015). {\it Bayesian Inference for Partially Identified Models: Exploring the Limits of Limited Data.}  CRC Press, New York.

\bibitem{} Haile, P., \& E. Tamer (2003). Inference with an incomplete model of English auctions. {\it Journal of Political Economy} 111, 1–51.  

\bibitem{} Kitagawa, T. (2012). Estimation and inference for set-identiﬁed parameters using posterior lower probability. {\it Working paper, University College London.}
http://www.homepages.ucl.ac.uk/~uctptk0/Research/LowerUpper.pdf

\bibitem{} Kline, B., \& Tamer, E. (2016). Bayesian inference in a class of partially identified models. {\it Quantitative Economics}. (To appear.)

%Another consequence is that, under certain necessary and sufficient conditions, the
%(1−α)-level credible set for the identified set is also an (1−α)-level frequentist confidence
%set for the identified set. This result means that there is an “asymptotic equivalence”
%between Bayesian and frequentist approaches to partially identified models, if the focus
%is on inference concerning the identified set rather than the partially identified parameter,
%which was the focus in other results including Moon and Schorfheide (2012).

\bibitem {} Liao, Y. (2010). {\it Bayesian Analysis in Partially Identified Parametric and Nonparametric Models}. Ph.D. thesis, Northwestern University.

%\bibitem{}  Liao, Y. \& W. Jiang (2010). Bayesian analysis in moment inequality models.
%{\it The Annals of Statistics} 38, 275--316.

\bibitem{} Manski, C. (2003). {\it Partial Identification of Probability Distributions}. Springer-Verlag, New York.
 

\bibitem{} Manski, C., \& E. Tamer (2002).   Inference on regressions with interval data on a regressor or outcome. {\it Econometrica} 70, 519–547.  

\bibitem{}  Moon, H. R. \&  F. Schorfheide (2012).  
Bayesian and frequentist inference in partially identified  models.
 \textit{Econometrica} 80,   755–782. 

\bibitem{}  Poirier, D. J. (1998). Revising beliefs in nonidentified models.
\textit{Econometric Theory} 14,  483-509. 

\bibitem{} van der Vaart, A. W. (2000). {\it Asymptotic statistics}. Cambridge University Press.

 

\end{thebibliography}
\end{document}